\documentstyle{amsart}

\def\N{{\Bbb N}}

\def\Q{{\Bbb Q}}
\def\R{{\Bbb R}}
\def\C{{\Bbb C}}
\def\P{{\Bbb P}}
\def\F{{\Bbb F}}
\def\G{{\Bbb G}}
\def\M{{\Bbb M}}

\def\A{{\cal A}}

\def\D{{\cal D}}
\def\H{{\cal H}}
\def\O{{\cal O}}
\def\T{{\cal T}}
\def\X{{\cal X}}
\def\calC{{\cal C}}
\def\calP{{\cal P}}
\def\calM{{\cal M}}

\def\x{{\vec{x}}}
\def\z{{\vec{z}}}
\def\g{{\frak{g}}}
\def\p{{\frak{p}}}

\def\Atilde{\widetilde{\A}}

\def\Btilde{\widetilde{B}}

\def\Ptilde{\widetilde{P}}
\def\Vtilde{\widetilde{V}}
\def\Xtilde{\widetilde{X}}
\def\Ytilde{\widetilde{Y}}

\def\Phat{{\widetilde{\cal P}}}
\def\Hcomp{{\widehat{\H}}}
\def\Otilde{{\widetilde{\O}}}
\def\Omegatilde{{\widetilde{\Omega}}}
\def\phibar{\overline{\phi}}
\def\Xtilde{\widetilde{X}}
\def\toptilde{\widetilde{\text{Top}}}

\def\dot{\bullet}
\def\vol{\mathrm{vol}}
\def\blank{\underline{\phantom{x}}}

\newcommand\spn{\operatorname{span}}
\newcommand\im{\operatorname{im}}
\newcommand\cone{\operatorname{Cone}}
\newcommand\Hom{\operatorname{Hom}}
\newcommand\Ext{\operatorname{Ext}}
\newcommand\comptensor{\operatorname{\widehat{\otimes}}}

\def\Hcts{H_{cts}}
\def\Hdel{H_{\D}}
\def\Hmd{H_{\calM\D}}

\def\Homcts{\Hom^{cts}}

\newtheorem{theorem}{Theorem}[section]
\newtheorem{lemma}[theorem]{Lemma}
\newtheorem{proposition}[theorem]{Proposition}
\newtheorem{corollary}[theorem]{Corollary}

\theoremstyle{definition}

\newtheorem{definition}[theorem]{Definition}
\newtheorem{example}[theorem]{Example}

\theoremstyle{remark}

\newtheorem{remark}[theorem]{Remark}

\begin{document}

\title{The Existence of Higher Logarithms}

\author{Richard M.~Hain}

\address{Department of Mathematics\\
Duke University\\ Durham, NC 27708-0320}

\thanks{Research  supported in part by grants from the National
Science Foundation}

\email{hain@@math.duke.edu}
\date{ August, 1993}

\maketitle

\begin{abstract}
In this paper we establish the existence of all higher logarithms
as Deligne cohomology classes in a sense slightly weaker than that of
\cite[\S 12]{hain-macp}, but in a sense that should be strong enough

for defining Chern classes on the algebraic $K$-theory of complex
algebraic varieties. In
particular, for each integer $p \ge 1$, we construct a multivalued
holomorphic function on a Zariski open  subset of the self dual
grassmannian of $p$-planes in $\C^{2p}$ which satisfies a
canonical $2p+1$ term functional equation. The key new technical
ingredient is the construction of a topology on the generic
part of each Grassmannian which is coarser than the Zariski
topology and where each open contains another which is both
a $K(\pi,1)$ and a rational $K(\pi,1)$.
\end{abstract}

\section{Introduction}\label{intro}

Denote by $G^p_q$ the Zariski open subset of the grassmannian of
$q$-dimensional linear subspaces of $\P^{p+q}$ which are transverse
to each coordinate hyperplane and each of their intersections.
Intersecting elements of $G^p_q$ with each of the $p+q+1$ coordinate
hyperplanes of $\P^{p+q}$ defines $p+q+1$ maps
$$
A_i : G^p_q \to G^p_{q-1}, \qquad 0\le i \le p+q.
$$
The spaces $G^p_q$ with $0 \le q \le p$ and the face maps $A_i$
form a truncated simplicial variety $G^p_\bullet$.

In \cite[\S 12]{hain-macp} (see also \cite{b-mcp-s}) the $p$th higher
logarithm is defined as a certain element of the
``multivalued Deligne cohomology'' of $G^p_\bullet$.
In that paper the existence of only the first three higher
logarithms was established.

In this paper we establish the existence of all higher logarithms,
but in a sense slightly weaker than that of \cite{hain-macp} ---

we show that for each $p$, there is a Zariski open subset
$U^p_\bullet$
of $G^p_\bullet$ on which the $p$th higher logarithm is defined as a
multivalued (and ordinary) Deligne cohomology class. This should be
sufficient to show that the $p$th Chern classes on the algebraic
$K$-theory of a complex algebraic variety is represented by the
$p$th higher logarithm (cf.\ \cite{goncharov:chern},
\cite{hain-yang}).
The key new technical
ingredient is the construction of a topology on the generic
part of each Grassmannian which is coarser than the Zariski
topology and where each open contains another which is both
a $K(\pi,1)$ and a rational $K(\pi,1)$.

Hanamura and MacPherson \cite{hanamura-macp_2} have a geometric
construction of the part of each higher logarithms that lies in
the multivalued de Rham complex of $G^p_\bullet$. The part of our
higher logarithm that lies in the multivalued de~Rham complex of
$G^p_\dot$ is only defined generically, so their result is stronger
than ours in this respect (cf. Remark \ref{hana-mac-const}),
but our result is stronger than theirs in that we construct higher
logarithms as both multivalued and ordinary Deligne cohomology
classes.

One part of the cocycle defining the $p$th higher logarithm is a
multivalued function $L_p$ defined on the Zariski open subset
$U^p_{p-1}$
of the self dual Grassmannian of $p$ planes in $\C^{2p}$.  The
cocycle
condition implies that this multivalued function satisfies the
canonical
$2p+1$ term functional equation
$$
\sum_{i=0}^{2p} (-1)^i A_i^\ast L_p = 0.
$$
In the cases $p=1,2,3$, the function
$L_p$ has a single valued cousin $D_p$ which also satisfies the
functional equation
$$
\sum_{i=0}^{2p} (-1)^i A_i^\ast D_i = 0.
$$
The first function $D_1$ is simply $\log|\phantom{x}|$, the second
is the Bloch-Wigner function, and the third, whose existence was
established in \cite[\S 11]{hain-macp}, can be expressed in terms of
Ramakrishnan's single valued cousin of the classical trilogarithm,
as was proved by Goncharov \cite{goncharov:trilog}. The functional
equation implies that $D_p$ ($p=1,2,3$) represents an element of
$H^{2p-1}(GL_p(\C),\C/\R(p))$.  This class is known to be a non-zero
rational multiple
of the Borel element, the class used to define the
Borel regulator (\cite{bloch}, \cite{dupont},
\cite{goncharov:trilog},
\cite{yang}, see also \cite{hain:poly}). The single valued cousins
of the higher logarithms constructed in this paper are constructed
in \cite{hain-yang} where it is shown that each represents
a non-zero rational multiple of the Borel class.

We now discuss the content of this paper in more detail.  As in
\cite{hain-macp}, the algebra of multivalued differential forms
on an algebraic manifold will be denoted by $\Omegatilde^\bullet(X)$.
We
will usually denote the ring of multivalued functions
$\Omegatilde^0(X)$
by $\Otilde(X)$. There is a weight filtration $W_\bullet$ on
$\Omegatilde^\bullet(X)$ which gives it the structure of a filtered
d.g.\ algebra. The category of complex algebraic manifolds and
regular
maps between them will be denoted by $\A$.  Since the pullback of a
multivalued function under a regular map $X\to Y$ is not well
defined, it is necessary to refine the category $\A$ in order that
$\Omegatilde^\bullet$ becomes a well defined functor. Such a
refinement $\Atilde$ of $\A$ is defined in \cite[\S 2]{hain-macp}.
The objects of $\Atilde$ are universal coverings $\widetilde{X}\to X$
of objects of $\A$, and the morphisms are commutative squares
$$
\begin{array}{ccc}
\widetilde{X} & \to & \widetilde{Y} \\
\downarrow & & \downarrow \\
X & \to & Y
\end{array}
$$
where the bottom arrow is a morphism of $\A$. The truncated
simplicial variety $G^p_\bullet$ has a natural lift to a
simplicial object of $\Atilde$, \cite[(5.4)]{hain-macp}.

Denote the Deligne-Beilinson cohomology of a smooth (simplicial)
variety $X$ by $\Hdel^\bullet(X,\Q(p))$. In \cite[\S 12]{hain-macp}
the multivalued Deligne cohomology of a simplicial object $X_\bullet$
of $\Atilde$ was defined. It will be denoted by
$\Hmd^\bullet(X_\bullet,\Q(p))$.

There are several equivalent ways to define rational $K(\pi,1)$
spaces,
but for our
purposes in the introduction, perhaps the most pertinent comment is
that a Zariski open subset $U$ of a grassmannian is a rational
$K(\pi,1)$ if and only if
$$
H^\bullet(W_l\Omegatilde^\bullet(U)) = \C
$$
for all $l\ge 0$. Using this, we show in Section \ref{pf_dunno} that
if $X_\bullet$ is a simplicial object of $\Atilde$ and each $X_q$ is
a rational $K(\pi,1)$, then there is a natural isomorphism
$$
\Hdel^\bullet(X_\bullet,\Q(p)) \approx \Hmd^\bullet(X_\bullet,\Q(p)).
$$
This was stated without proof in \cite[(12.3)]{hain-macp}.

The main idea of this paper is
to exploit this fact by replacing $G^p_\bullet$ by a Zariski open
subset $U^p_\bullet$ where each $U^p_q$ is a rational $K(\pi,1)$.
Once one has done this and established that $U^p_\dot$ lifts to a
simplicial object of $\Atilde$, the proof of the existence of higher
logarithms is relatively straightforward --- there is a natural
$GL_p(\C)$ bundle over $U^p_\bullet$ whose $p$th Chern class is an
element of
$$
\Hdel^{2p}(U^p_\bullet,\Q(p)) \approx \Hmd^{2p}(U^p_\bullet,\Q(p)).
$$
The $p$th higher logarithm is a suitable rational multiple of this
class.%
\footnote{If one only wants the multivalued function, or the higher
logarithm in the sense of \cite[(6.1)]{hain-macp}, then one can
appeal
directly to the analogue of \cite[(8.9)]{hain-macp} for $U^p_\bullet$
--- cf. (\ref{exist_1}).}

To recapitulate, one of the main obstacles to proving the existence
of the
$p$th higher logarithm is to establish the existence of such a
Zariski
open subset $U^p_\bullet$ of $G^p_\dot$ where each $U^p_q$ is a
rational
$K(\pi,1)$. If the Zariski topology had the property that each open
set
contains another open which is a rational $K(\pi,1)$, then one could
easily find the sought after open subset $U^p_\dot$ of $G^p_\dot$.
Unfortunately, this is not the case (cf. \cite[(9.7)]{hain:cycles}).
For this reason we introduce a coarser topology on the $G^p_q$,
called
the {\it constructible topology}, which does enjoy this property.
This
is done in Section \ref{topology}.

We conclude the introduction with a brief description of the
constructible
topology and the idea behind the proof of the existence of
$U^p_\dot$. The
first point is that each $\xi \in G^p_q$ determines an ordered
configuration
of $p+q+1$ points in $\P^{p-1}$, no $p$ of which lie on a hyperplane;
the configuration is well defined up to projective equivalence. To
see
how this works, note that the set of $(q+1)$-dimensional planes in
$\P^{p+q}$ which contain $\xi$ comprise a projective space of
dimension
$p-1$. The $j$th point of the configuration is the point of this
projective space which corresponds to the join of the $j$th standard
basis vector with $\xi$. Each such configuration determines a
configuration of hyperplanes in $\P^{p-1}$ --- the hyperplanes are
those
determined by the $(p-1)$-element subsets of the points. The
configuration
corresponding to an element of $G^3_2$ and the corresponding
arrangement
of hyperplanes in $\P^2$ are illustrated in Figure \ref{fiber}.

A configuration of hyperplanes in $\P^{p-1}$ corresponds to a central
configuration of hyperplanes in $\C^p$.
Denote the central configuration in $\C^p$ which corresponds to
$\xi \in G^p_q$ by $\calC(\xi)$.

One's natural instinct when trying to understand the topology of the
$G^p_q$ is to use the face maps $A_i : G^p_q \to G^p_{q-1}$ to study
them inductively. The ``standard mistake'' is to believe that all
such
face maps are fibrations. If they were, life would be easier, but
less interesting. It is worthwhile to see how the face maps fail to
be fibrations as it is relevant to the proof of the existence of

$U^p_\dot$. Observe that the fiber of the face map
$A_i : G^p_q \to G^p_{q-1}$ over the point $\xi \in G^p_{q-1}$ is the
complement of the arrangement $\calC(\xi)$ in $\C^p$.

The simplest example where a face map is not a fibration is
provided by any of the face maps $A_i : G^3_3 \to G^3_2$. The
projectivization of the generic fiber is the complement of an
arrangement
determined by six points in $\P^2$, no three of which lie on a line,
and
where no three of the lines they determine are concurrent, except at
one
of the points $x_0,\dots, x_5$. The complement of the arrangement on
the
right hand side of Figure \ref{fiber} is the projectivization of a
special fiber of $A_6 : G^3_3 \to G^3_2$ as there is an exceptional
triple
point. Since the topology of the fiber of $A_6 : G^3_3 \to G^3_2$ is
not
constant, $A_6$ is not a fibration.

\begin{figure}

\begin{picture}(400,200)
\put(10,0){\begin{picture}(200,200)
\put(80,150){\circle*{6}}		
\put(88,147){$x_0$}
\put(30,50){\circle*{6}}		
\put(10,57){$x_1$}
\put(130,50){\circle*{6}}		
\put(120,32){$x_2$}
\put(80,100){\circle*{6}}		
\put(65,77){$x_3$}
\put(180,100){\circle*{6}}		
\put(172,110){$x_4$}
\put(30,130){\circle*{6}}		
\put(15,114){$x_5$}
\end{picture}}

\put(250,0){\begin{picture}(200,200)
\put(80,150){\circle*{6}}		
\put(88,147){$x_0$}
\put(30,50){\circle*{6}}		
\put(10,57){$x_1$}
\put(130,50){\circle*{6}}		
\put(120,32){$x_2$}
\put(80,100){\circle*{6}}		
\put(65,77){$x_3$}
\put(180,100){\circle*{6}}		
\put(172,110){$x_4$}
\put(30,130){\circle*{6}}		
\put(15,114){$x_5$}

\put(10,10){\line(1,2){90}}		
\put(80,-10){\line(0,1){200}}		
\put(150,10){\line(-1,2){90}}		
\put(10,30){\line(1,1){160}}		
\put(10,50){\line(1,0){180}}		
\put(170,10){\line(-1,1){160}}		

\put(200,100){\line(-1,0){190}}		
\put(200,90){\line(-2,1){190}}		
\put(201,107){\line(-3,-1){195}}	
\put(200,120){\line(-1,-1){130}}	

\put(10,122){\line(5,2){170}}		
\put(30,10){\line(0,1){180}}		
\put(10,146){\line(5,-4){170}}		
\put(10,142){\line(5,-3){180}}		
\put(200,96){\line(-5,1){190}}		
\end{picture}}
\end{picture}

\caption{}\label{fiber}
\end{figure}

The basic closed subsets of $G^p_q$ in the constructible topology
are defined to be the closure of the set of points $\xi$ where the
combinatorics of $\calC(\xi)$ is fixed.
For example, the closure of the set of points in $G^3_2$ where the
lines
$x_0x_2$, $x_1x_3$ and $x_4x_5$ intersect in a single point
(as in Figure \ref{fiber}) is a closed
subset of $G^3_2$ in the constructible topology. The combinatorial
objects which parameterize the closed sets are called {\it
templates}.

Observe that $A_0 : G^3_3 \to G^3_2$ is a fibration over the
constructible open subset of $G^3_2$ which consists of all $\xi$
for which the projectivization of $\calC(\xi)$ contains no
exceptional
triple points. By passing to a constructible open subset of $G^3_3$,
one can arrange for the generic fiber of $A_0$ to be the complement
of an arrangement of fiber type, and by restricting $A_0$ to a
smaller
constructible open subset of $G^3_2$ we may assume that
$A_0$ is a fibration whose fiber is the complement of an arrangement
of fiber type. Since arrangements of fiber type are rational
$K(\pi,1)$s, this provides, via (\ref{fibration}), the inductive step
needed for finding the open subset $U^p_\dot$ of $G^p_\dot$ in which
each $U^p_q$ is a rational $K(\pi,1)$.

It is assumed that the reader is familiar with \cite{hain-macp}.
\medskip

\noindent{\sl Conventions.} In this paper, all simplicial objects are
strict --- that is, they are functors from the category $\Delta$ of
finite ordinals and {\it strictly} order preserving maps to, say,
the category of algebraic varieties.

As is standard, the finite set $\{0,1,\dots,n\}$ with its natural
ordering will be denoted by $[n]$. Let $r$ and $s$ be positive
integers
with $r \le s$. Denote the full subcategory of $\Delta$ whose objects
are the ordinals $[n]$ with $r\le n \le s$ by $\Delta[r,s]$. An
$(r,s)$
{\it truncated} simplicial object of a category $\calC$ is a
contravariant
functor from $\Delta[r,s]$ to $\calC$.

The word {\it simplicial} will be used generically to refer to both
simplicial objects and truncated simplicial objects.

By Deligne cohomology, we shall mean Beilinson's refined version
of Deligne cohomology as defined in \cite{beilinson:ahc}. It can
be expressed as an extension
$$
0 \to \Ext^1_\H(\Q,H^{k-1}(X,\Q(p))) \to \Hdel^k(X,\Q(p)) \to
\Hom_\H(\Q,H^k(X,\Q(p)))\to 0
$$
where $\H$ denotes the category of $\Q$ mixed Hodge structures.
\medskip

\noindent{\sl Acknowledgements.} I would like to thank Jun Yang for
helpful discussions regarding material in Sections \ref{existence}
and
\ref{proof}.
\medskip

\section{Constructible Configurations and Templates}

Fix a ground field $\F$. Denote the projective space $\P^m(\F)$ over
$\F$  by $\P^m$. By a {\it configuration} of $n$ points in $\P^m$, we
shall mean an  element $\x$ of $\left(\P^m\right)^n$. A
{\it subconfiguration} of $\x$ is any element of $(\P^m)^l$, $l\le
n$,
obtained by deleting some of the components of $\x$.

A {\it  linear configuration} in $\P^m$ is a finite collection of
linear subspaces of $\P^m$.  The {\it complete configuration}
$\Hcomp$ associated to a linear configuration $\H = \{L_1,\ldots, L_l
\}$
in $\P^m$ is the configuration consisting of the $L_j$
and all of their non-empty intersections. A linear
configuration is {\it complete} if it equals its completion.

The {\it join} of two linear subspaces $L_1$, $L_2$ of $\P^m$ is the
smallest linear subspace of $\P^m$ which contains them both. It will
be denoted by $L_1 \ast L_2$.

\begin{definition}\label{derived}
The set $\D(\x)$ of linear configurations in $\P^m$
{\it derived} from a particular configuration $\x$ of $n$ points
in $\P^m$ is the unique set of linear configurations in $\P^m$
which satisfies the following properties:
\begin{enumerate}
\item  the completion of the configuration consisting of all
hyperplanes
in $\P^m$ that are spanned by a subconfiguration of $\x$ is in
$\D(\x)$;
\item every $\H \in \D(\x)$ is complete;
\item if $\H \in \D(\x)$,  $L\in \H$ and $\X$ is a subconfiguration
of $\x$ such that $L \ast \spn \X$ is a hyperplane, then the
completion of $\H \cup \{L \ast \spn \X\}$ is also in $\D(\x)$.
\end{enumerate}
\end{definition}

\begin{figure}[htp]

\begin{picture}(400,200)
\put(10,0){\begin{picture}(200,200)
\put(80,150){\circle*{6}}		
\put(75,157){$x_0$}
\put(30,50){\circle*{6}}		
\put(25,57){$x_1$}
\put(130,50){\circle*{6}}		
\put(125,57){$x_2$}
\put(80,100){\circle*{6}}		
\put(75,107){$x_3$}
\put(180,100){\circle*{6}}		
\put(175,107){$x_4$}
\end{picture}}

\put(250,0){\begin{picture}(200,200)
\put(80,150){\circle*{6}}		
\put(30,50){\circle*{6}}		
\put(130,50){\circle*{6}}		
\put(80,100){\circle*{6}}		
\put(180,100){\circle*{6}}		

\put(10,10){\line(1,2){90}}		
\put(80,-10){\line(0,1){200}}		
\put(150,10){\line(-1,2){90}}		
\put(10,30){\line(1,1){120}}		
\put(10,50){\line(1,0){180}}		
\put(170,10){\line(-1,1){140}}		

\put(200,100){\line(-1,0){180}}		
\put(200,90){\line(-2,1){150}}		
\put(201,107){\line(-3,-1){195}}	
\put(200,120){\line(-1,-1){130}}	
\end{picture}}
\end{picture}

\caption{}\label{points}
\end{figure}

\begin{example}
Let $\x$ be the configuration $(x_0,x_1,x_2,x_3,x_4)$ of 5 points
in $\P^2(\R)$ depicted in the left half of Figure \ref{points}. The
right
half of Figure \ref{points} depicts the configuration defined
in (1) of the definition of $\D(\x)$.  Every other linear
configuration
$\H$ in $\D(\x)$ contains this configuration. The first linear
configuration depicted in Figure \ref{config} is in $\D(\x)$, the
second is not.
\end{example}

\begin{figure}

\begin{picture}(400,200)
\put(10,0){\begin{picture}(200,200)
\put(80,150){\circle*{6}}		
\put(30,50){\circle*{6}}		
\put(130,50){\circle*{6}}		
\put(80,100){\circle*{6}}		
\put(180,100){\circle*{6}}		

\put(10,10){\line(1,2){90}}		
\put(80,-10){\line(0,1){200}}		
\put(150,10){\line(-1,2){90}}		
\put(10,30){\line(1,1){120}}		
\put(10,50){\line(1,0){180}}		
\put(170,10){\line(-1,1){140}}		

\put(200,100){\line(-1,0){180}}		
\put(200,90){\line(-2,1){150}}		
\put(201,107){\line(-3,-1){195}}	
\put(200,120){\line(-1,-1){130}}	
\put(200,96){\line(-5,1){170}}		
\end{picture}}

\put(250,0){\begin{picture}(200,200)
\put(80,150){\circle*{6}}		
\put(30,50){\circle*{6}}		
\put(130,50){\circle*{6}}		
\put(80,100){\circle*{6}}		
\put(180,100){\circle*{6}}		

\put(10,10){\line(1,2){90}}		
\put(80,-10){\line(0,1){200}}		
\put(150,10){\line(-1,2){90}}		
\put(10,30){\line(1,1){120}}		
\put(10,50){\line(1,0){180}}		
\put(170,10){\line(-1,1){140}}		

\put(200,100){\line(-1,0){180}}		
\put(200,90){\line(-2,1){150}}		
\put(201,107){\line(-3,-1){195}}	
\put(200,120){\line(-1,-1){130}}	
\put (20,117){\line(1,0){190}}		
\end{picture}}

\end{picture}
\caption{}\label{config}
\end{figure}

The class of all order preserving functions $r : P \to \N$
from a partially ordered set into $\N$ forms a category $\calP$. A
morphism
$\Phi$ from $r_1 : P_1 \to \N$ to $r_2: P_2 \to \N$ is  an order
preserving function $\phi : P_1 \to P_2$ such that the diagram
$$
\begin{matrix}
P_1 & \stackrel{\phi}{\longrightarrow} & P_2 \\
r_1 \searrow & & \swarrow r_2 \\
& \N &
\end{matrix}
$$
commutes. We shall denote the {\it $k$-dimensional elements}
$r^{-1}(k)$
of $P$ by $P_k$.

\begin{definition} A {\it template} is an isomorphism class of
objects of the category $\calP$.
\end{definition}

Templates are a generalization of matroids.

Each linear configuration $\H$ in $\P^m$ is a partially ordered
set---the linear subspaces are ordered by inclusion.  Define
a rank function $r: \H \to \N$ by defining $r(L) = \dim L$ for each
$L \in \H$.
In this way we associate a template to each linear configuration.

The template of a linear configuration $\H$ derived from a
configuration of points $\x$ in $\P^m$ has additional structure;
namely, the marking of the points of $\x$.  For this reason, we
now define marked templates.

Denote the set $\{0,1,\ldots,n\}$ by $[n]$. One can consider the
class of triples $(P,r,\psi)$, where $P$ is a partially ordered
set, $r:P \to \N$ is an order preserving function, and where
$\psi : [n] \to P_0$ is a function.  These form a category $\Phat_n$;
the morphisms are order preserving maps which preserve the rank
functions $r$ and the markings $\psi$.

\begin{definition}
An {\it $n$-marked template} is an isomorphism class of objects
of the category $\Phat_n$.
\end{definition}

Marked templates are a generalization of oriented matroids.

Each linear configuration derived from a configuration $\x$ of
$n+1$ points in $\P^m$ determines an $n$-marked template.  If
$\x = (x_0,x_1,\ldots , x_n)$, then the  marking
$\psi : [n] \to \H_0$ is defined by $\psi(j) = x_j$. We will
view $\D(\x)$ as a set of
marked linear configurations. Denote the set of $n$-marked
templates
$$
\left\{ T(\H) : \H \in \D(\x)\right\}
$$
associated to a configuration $\x$ of $n+1$ points in $\P^m$
by $\T(\x)$.  Taking $\H$ to $T(\H)$ defines a bijection
$$
\D(\x) \to \T(\x).
$$
We shall denote the element of $\D(\x)$ which corresponds to
$T\in \T(\x)$ by $\H_T$.

The group of projective equivalences $PGL_{m+1}$ acts on the
set of linear configurations in $\P^m$. Observe that if two
linear configurations are projectively equivalent, they determine
the same template. Consequently, $\T(\x)$ depends only on the
projective equivalence class of $\x$.

\section{Hyperplane Arrangements of Fiber Type}

We retain the notation of the previous section. We inductively define
what it means for an arrangement of hyperplanes in $\F^n$ to be
of {\it fiber type}. First, every arrangement of distinct points
in $\F$ is of fiber type. An arrangement of hyperplanes $\H$ in
$\F^n$
is of fiber type if there is a linear projection $\phi:\F^n\to
\F^{n-1}$
and an arrangement of hyperplanes $\A$ in $\F^{n-1}$ of fiber type
such that
\begin{enumerate}
\item[(a)] the arrangement $\phi^{-1}\A$ is a sub-arrangement of
$\H$;
\item[(b)] the image under $\phi$ of each element of
$\H - \phi^{-1}\A$
is all of $\F^{n-1}$;
\item[(c)] for each $u\in \F^{n-1}-\cup\A$, the number of points in
the induced arrangement of points $\phi^{-1}(u)\cap \H$ of
$\phi^{-1}(u)$ by $\H$ is independent of $u$.
\end{enumerate}

When $\F$ is $\R$ or $\C$, the conditions (a), (b) and (c) imply that
the projection $\psi : \C^n - \bigcup\H \to \C^{n-1} - \bigcup \A$ is
a topological fiber bundle.

\begin{proposition}\label{fiber_type}
If $\x$ is a configuration of $n$ points in $\P^m$, then, for
each $T \in \T(\x)$, there is $T'\in \T(\x)$ such that
$\H_T \subseteq \H_{T'}$ and such that $\H_{T'}$ is an arrangement
of fiber type.
\end{proposition}

\begin{pf}
We prove the result by induction on $m$. The result is trivially true
when $m=1$.  Now suppose that $m\ge 1$.  The image of
$(x_0,\ldots, x_{n-1})$ under the linear projection
$$
\phi : \P^m - \{x_n\} \to \P(T_{x_n}\P^m) \approx \P^{m-1}
$$
is an $(n-1)$-marked configuration $\z$ of points in $\P^{m-1}$.
The linear subspaces $L\in \H_T$ which contain $x_0$ induce a linear
arrangement $\overline{\H}_T$ in $\P(T_{x_n}\P^m)$. It is easy to
check that $\overline{\H}_T \in \D(\z)$. By induction, there is a
template $S\in \T(\z)$ such that the arrangement $\H_S$ in
$\P(T_{x_n}\P^m)$ is of fiber type and contains $\overline{\H}_T$.
The inverse image of $\H_S$ under
$\phi$ is an arrangement $\widetilde{\H}_S$ of hyperplanes in
$\P^m$ each of whose hyperplanes contains $x_n$. The completion
of the linear arrangement
$$
\H := \widetilde{\H}_S \cup \H_T
$$
is an element of $\D(\x)$. The projection $\phi$ induces a linear
projection
$$
\psi : \P^m - \cup \H \to \P(T_{x_n}\P^m) - \cup \H_S
$$
whose fibers are punctured lines. Adding to $\H$ the hyperplanes
in $\P^m$ which are the join of $x_n$ with a codimension 2 stratum
of $\H$ we obtain a linear arrangement $\H'$ in $\P^m$ such that
the restriction of $\psi$ to $\P^m - \cup \H'$ is a linear
map
$$
\P^m - \cup \H' \to \P(T_{x_n}\P^m) - \cup \H_S
$$
each of whose fibers is $\P^1$ minus the same number of points. That
is,
the arrangement $\H'$ is of fiber type. Let $T'\in \T(\x)$ be the
template which corresponds to the completion of $\H'$. Then
$\cup \H_{T'} = \cup \H'$, and so $\H_{T'}$ is an arrangement of
fiber type which contains $\H_T$.
\end{pf}

\section{The Generic Grassmannians}\label{topology}

As in the previous sections, $\F$ will denote a fixed ground field.
First recall that the grassmannian $G(q,\P^{p+q})$ of
$q$-dimensional subspaces of $\P^{p+q}$ can be viewed as the
orbit space
$$
\left\{ (v_0,v_1,\ldots,v_{p+q}) \in (\F^p)^{p+q+1} :
v_0, \ldots, v_{p+q} \text{ span } \F^p\right\}/GL_p(\F)
$$
where $GL_p$ acts diagonally (cf. \cite[\S 5]{hain-macp}).

The generic part $G^p_q$ of $G(q,\P^{p+q})$ is defined to be the set
of
those points in $G(q,\P^{p+q})$ which correspond to $(p+q+1)$-tuples
of vectors $(v_0,\ldots,v_{p+q})$ in $\F^p$ where each $p$ of the
vectors are linearly independent.

The torus
$$
(\G_m)^{p+q} \approx (\G_m)^{p+q+1}/\text{ diagonal}
$$
acts on $G^p_q$ via the action
$$
(\lambda_0,\ldots,\lambda_{p+q}) : (v_0,\ldots,v_{p+q}) \mapsto
(\lambda_0 v_0,\ldots,\lambda_{p+q} v_{p+q}).
$$
The quotient space is the variety
$$
Y^p_q := \left\{(x_0,\ldots,x_{p+q}) \in (\P^{p-1})^{p+q+1} :
\text{each $p$ of the points span $\P^{p-1}$} \right\}/PGL_p.
$$
The morphism $\pi : G^p_q \to Y^p_q$ is a principal
$(\G_m)^{p+q}$-bundle with a section \cite[(5.9)]{hain-macp}.
Consequently,
$$
G^p_q \approx Y^p_q \times (\G_m)^{p+q}.
$$

Denote the point of $\P(V)$ which corresponds to $v\in V-\{0\}$
by $[v]$. The point $v$ of $G^p_q$ corresponding to the orbit of
$(v_0,v_1,\ldots,v_{p+q})$ determines the point
$$
\x(v) = ([v_0],[v_1],\ldots,[v_{p+q}])
$$
of $Y^p_q$.  We can therefore associate to each
point of $G^p_q$ the set $\T(\x(v))$ of $(p+q)$-marked templates.

For each $(p+q)$-marked template $T$, define the subset $E^p_q(T)$
of $G^p_q$ to be the Zariski closure of
$$
\left\{v \in G^p_q : T \in \T(\x(v)) \right\}.
$$
We will define two templates $T_1$  and $T_2$ to
be {\it $(p,q)$-equivalent} if the subvarieties $E^p_q(T_1)$ and
$E^p_q(T_2)$ of $G^p_q$ are equal.

\begin{example}
The two configurations in Figure \ref{templates} determine
templates $T_1$ and $T_2$, respectively.  Both $E^3_4(T_1)$ and
$E^3_4(T_2)$ are  proper subvarieties of $G^3_4$, and $E^3_4(T_2)$
is a proper subvariety of $E^3_4(T_1)$.
\end{example}

\begin{figure}[htp]
\begin{picture}(440,220)(0,0)
\put(20,0){\begin{picture}(220,220)
\put(50,50){\circle*{6}}
\put(50,150){\circle*{6}}
\put(150,150){\circle*{6}}
\put(150,50){\circle*{6}}
\put(20,125){\circle*{6}}
\put(180,125){\circle*{6}}
\put(100,35){\circle*{6}}
\put(100,165){\circle*{6}}

\put(20,20){\line(1,1){160}}
\put(0,125){\line(1,0){200}}
\put(20,180){\line(1,-1){160}}
\put(100,10){\line(0,1){180}}
\end{picture}}

\put(240,0){\begin{picture}(220,220)
\put(50,50){\circle*{6}}
\put(50,150){\circle*{6}}
\put(150,150){\circle*{6}}
\put(150,50){\circle*{6}}
\put(20,100){\circle*{6}}
\put(180,100){\circle*{6}}
\put(100,35){\circle*{6}}
\put(100,165){\circle*{6}}

\put(20,20){\line(1,1){160}}
\put(0,100){\line(1,0){200}}
\put(20,180){\line(1,-1){160}}
\put(100,10){\line(0,1){180}}
\end{picture}}
\end{picture}
\caption{}\label{templates}
\end{figure}

Define $F^p_q(T) \subseteq Y^p_q$ to be the quotient of $E^p_q(T)$
by the torus action. Observe that:

\begin{proposition}\label{prod}
For each template $T$, the varieties
$E^p_q(T)$ and $F^p_q(T) \times (\G_m)^{p+q}$ are isomorphic. \qed
\end{proposition}

\begin{definition}
The {\it constructible topology} on $G^p_q$ is the topology whose
closed sets are finite unions of the sets $E^p_q(T)$. The
{\it constructible topology} on $Y^p_q$ is the topology whose closed
sets are finite unions of the sets $F^p_q(T)$. The constructible
topology on a subset of $G^p_q$ or $Y^p_q$ is the topology
induced from the constructible topology on $G^p_q$ or $Y^p_q$.
In particular, the sets $E^p_q(T)$ and $F^p_q(T)$ have constructible
topologies.
\end{definition}

Evidently, the closed subsets of $G^p_q$ are precisely the inverse
images of closed subsets of $Y^p_q$ under the projection
$G^p_q \to Y^p_q$. Note that the constructible topology is
coarser than the Zariski topology.

\begin{proposition}\label{generic}
For each $(p,q)$-marked template $T$,
$$
\left\{v \in F^p_q(T) : T \in \T(\x(v))\right\}
$$
is a constructible open subset of $F^p_q(T)$.  \qed
\end{proposition}

\section{Rational $K(\pi,1)$ Spaces}\label{rat_space}

In this section we briefly review the definition and basic properties
of rational $K(\pi,1)$ spaces.  Relevant references include
\cite{kohno}, \cite{falk}, \cite{hain:cycles} and \cite{hain-macp}.

As motivation, recall that if a topological space $X$ is a
$K(\pi,1)$,
then there is a natural isomorphism
$$
H^\bullet(X,\M) \approx H^\bullet(\pi_1(M),M)
$$
where $M$ is a $\pi_1(M,\ast)$ module, and $\M$ denotes the
corresponding
local system over $X$.

One can define the continuous cohomology of a group
$\pi$ by
$$
\Hcts^\bullet(\pi;\Q) =
\lim_{\stackrel{\to}{\Gamma}} H^\bullet(\Gamma,\Q)
$$
where $\Gamma$ ranges over all finitely generated nilpotent quotients
of $\pi$. There is an evident map
$$
\Hcts^\bullet(\pi,\Q) \to H^\dot(\pi,\Q).
$$
A topological space $X$ is defined to be a rational $K(\pi,1)$ if
the composition
$$
\Hcts^\dot(\pi_1(X),\Q) \to H^\dot(\pi_1(X),\Q) \to H^\dot(X,\Q)
$$
is an isomorphism. Every nilmanifold is clearly both a $K(\pi,1)$
and a rational $K(\pi,1)$. In particular, the circle is both
a $K(\pi,1)$ and a rational $K(\pi,1)$.

The following results will be used in Section \ref{main_thm}.
Proofs of them can be found in \cite[\S5]{hain:cycles}.

\begin{theorem}\label{wedge}
The one point union of two rational $K(\pi,1)$s is a rational
$K(\pi,1)$. In particular, every Zariski open subset of $\C$ is
both a $K(\pi,1)$ and a rational $K(\pi,1)$. \qed
\end{theorem}

\begin{theorem}\label{fibration}
Suppose that $f:X \to Y$ is a fiber bundle with fiber $F$. If $Y$ and
$F$ are rational $K(\pi,1)$s, and if the natural action of
$\pi_1(Y,y)$ on each cohomology group of $F$ is unipotent, then
$X$ is a rational $K(\pi,1)$. \qed
\end{theorem}

Since each Zariski open subset of $\C$ is both a $K(\pi,1)$ and
a rational $K(\pi,1)$, and since the monodromy representations
associated to a linear fibration is trivial
\cite[(5.12)]{hain:cycles}
we obtain the following result:

\begin{corollary}\label{fiber-type}
The complement of an arrangement of hyperplanes in $\C^n$
which is of fiber type is both a $K(\pi,1)$ and a rational
$K(\pi,1)$.
\qed
\end{corollary}

\section{The Main Theorem}\label{main_thm}

In this section, we prove the following result.

\begin{theorem}\label{main}
Each constructible open subset of $G^p_q(\C)$ contains a
constructible
open subset which is both a $K(\pi,1)$ and a rational $K(\pi,1)$.
\end{theorem}

\begin{remark}
It is easy to show that in the cases of $G^p_0$ and $G^p_1$, the
constructible topology is trivial. That
is, the only constructible open sets in these spaces are the the
empty set and the whole space. Thus Theorem \ref{main} implies
that $G^p_0(\C)$ and $G^p_1(\C)$ are $K(\pi,1)$s and rational
$K(\pi,1)$s. This is clear in the case of $G^p_0$, and is proved
directly in the case of $G^p_1$ in \cite[(8.6)]{hain-macp}.
\end{remark}

The proof of Theorem \ref{main} occupies the rest of this section.
Since the classes of rational $K(\pi,1)$s and $K(\pi,1)$s are
closed under products, and since each constructible open subset
of $G^p_q$ is a product of the corresponding constructible open
subset
of $Y^p_q$ with $(\C^\ast)^{p+q}$, we need only prove that each
constructible open subset of $Y^p_q$ contains a constructible open
set
which is a both a $K(\pi,1)$ and a rational $K(\pi,1)$.

Suppose that $0\le i \le {p+q}$.  The {\it $i$th face map},
$$
A_i : Y^p_q \to Y^p_{q-1}
$$
is defined by forgetting the $i$th point of a configuration of
$p+q$ points in $\P^{p-1}$. The $i$th {\it dual face map}
$$
B_i : Y^p_q \to Y^{p-1}_q
$$
is obtained by projecting all but the $i$th point of a configuration
of
$p+q$ points in $\P^{p-1}$ onto a generic $\P^{p-2}$ using the $i$th
point as the center of the projection.

\begin{proposition}\label{cont}
If $T$ is a $(p+q)$-marked template, then for each integer $i$
satisfying
$0\le i \le p+q+1$, there is a $(p+q+1)$-marked template $A^iT$
(resp.\
$B^iT$) whose $(p,q+1)$-equivalence class (resp.\
$(p+1,q)$-equivalence
class) depends only on the $(p,q)$-equivalence class of $T$.
Moreover,
$$
A_i^{-1}F^p_q(T) = F^p_{q+1}(A^iT),\quad
A_i^{-1}E^p_q(T) = E^p_{q+1}(A^iT)
$$
and
$$
B_i^{-1}F^p_q(T) = F^{p+1}_{q}(B^iT),\quad
B^iE^p_q(T) = E^{p+1}_{q}(B^iT).
$$
In particular, the face maps and dual face maps are continuous
with respect to the constructible topology.
\end{proposition}

\begin{pf}
For simplicity of notation, we take $i=p+q+1$.  Suppose that $T$ is a
$(p+q)$-marked template. Represent it by the object $(P,r,\psi)$ of
the category $\Phat_{p+q}$. Denote the marked elements
$\psi(0), \ldots,\psi(p+q)$ of $P_0$ by $p_0,\ldots,p_{p+q}$.
For each subset $I$ of $\{0,\ldots,p+q\}$, denote the element
of $P$ which is the least upper bound of $\{p_i: i\in I\}$ by
$p_I$. Set $r_I = r(p_I)$. Define $A^i T$ and $B^i T$ both to
be the isomorphism class of the completion of the marked ordered
set obtained from $(P,r,\psi)$ by adding one extra element
$p_{p+q+1}$ to $P_0$, and elements $p_I\ast p_{p+q+1}$ to
$P_{1+r_I}$. This is the ``smallest'' template $T'$ for
which $A_iT' = T$.

For example, if $T$ is the 6-marked template associated
to the configuration on the left hand side of Figure \ref{invim},
then $A^6T$ is the 7-marked template associated to the
configuration on the right hand side of Figure \ref{invim}.

\begin{figure}
\begin{picture}(440,220)(0,0)
\put(20,0){\begin{picture}(220,220)
\put(50,50){\circle*{6}}
\put(50,150){\circle*{6}}
\put(150,150){\circle*{6}}
\put(150,50){\circle*{6}}
\put(100,40){\circle*{6}}
\put(100,160){\circle*{6}}

\put(20,20){\line(1,1){160}}
\put(20,180){\line(1,-1){160}}
\put(100,-10){\line(0,1){220}}
\end{picture}}

\put(280,0){\begin{picture}(220,220)
\put(50,50){\circle*{6}}
\put(50,150){\circle*{6}}
\put(150,150){\circle*{6}}
\put(150,50){\circle*{6}}
\put(100,40){\circle*{6}}
\put(100,160){\circle*{6}}

\put(0,100){\circle*{6}}
\put(-6,112){$x_6$}

\put(20,20){\line(1,1){160}}
\put(20,180){\line(1,-1){160}}
\put(100,-10){\line(0,1){220}}

\put(-20,80){\line(1,1){130}}
\put(-20,88){\line(5,3){170}}
\put(-21,93){\line(3,1){200}}

\put(-20,120){\line(1,-1){130}}
\put(-20,112){\line(5,-3){170}}
\put(-21,107){\line(3,-1){200}}
\end{picture}}
\end{picture}
\caption{}\label{invim}
\end{figure}

For $v\in Y^{p+1}_q$, it is clear that $T \in \T(\x(A_iv))$ if
and only if $A^i T \in \T(\x(v))$. It follows that
$$
A_i^{-1}F^p_q(T) = F^{p+1}_q(A^iT)
$$
and that the $(p+1,q)$-equivalence class of $A^i T$ depends
only on the $(p,q)$-equivalence class of $T$.

The corresponding statement for $E^p_q(T)$ follows from
(\ref{prod}). The statements with $A$ replaced by $B$
follow using the dual argument.
\end{pf}

The following definition is an analogue of (\ref{derived}) for
templates. It is used only in the proof of the next result.

\begin{definition}\label{generate}
Suppose that $(P,r,\psi)$ is an object of the category
$\Phat_n$.
Suppose that $f: [m] \to [n]$ is an order preserving
injection. Define the subset $Q$ of $P$ {\it generated} by $f$ to be
the smallest subset of $P$ which contains
$\left\{\psi\circ f(j) : j\in [m]\right\}$ and is
closed under the following operations:
\begin{enumerate}
\item if $S\subseteq \im f$, then the least upper bound of $S$ in
$P$ is in $Q$;
\item if $S\subseteq Q$, then the greatest lower bound of $S$ in $P$
is in $Q$;
\item if $v\in Q$ and $S\subseteq \im f$, then the greatest lower
bound of $S$ and $v$ in $P$ is in $Q$.
\end{enumerate}
\end{definition}

\begin{proposition}\label{closed}
If $T$ is a $(p+q)$-marked template, then for each integer $i$
satisfying
$0\le i \le p+q$, there is a $(p+q-1)$-marked template $A_iT$ (resp.\
$B_iT$) whose $(p,q-1)$-equivalence class (resp.\
$(p-1,q)$-equivalence
class) depends only on the $(p,q)$-equivalence class of $T$.
Moreover,
$A_iF^p_q(T)$ is a constructible open subset of
$F^p_{q-1}(A_iT)$, $A_iE^p_q(T)$ is a constructible open subset of
$E^p_{q-1}(A_iT)$, $B_iF^p_q(T)$ is a constructible open subset of
$F^{p-1}_{q}(B_iT)$ and $B_iE^p_q(T)$ is a constructible open
subset of $E^{p-1}_{q}(B_iT)$.
\end{proposition}

\begin{pf}
Suppose that $T$ is a $(p+q)$-marked template.  Let $(P,r,\psi)$ be
an object of the category $\Phat_{p+q}$ which represents $T$. Let
$Q$ be the partially ordered subset of $P$ generated by the $i$th
face map $d_i : [p+q] \to [p+q+1]$---that is, the unique order
preserving
 injection which omits the value $i$.  Let $(Q,r,\psi\circ d_i)$ be
the
object of $\Phat_{p+q-1}$ where $r$ is the restriction of the
rank function of $P$. Define $A_{p+q}T$ to be the
$(p+q-1)$-marked template which is represented by $(Q,r,\psi)$.

For example,  if $T$ is the template associated to the configuration
in the left hand side of Figure \ref{less}, then $A_2 T$ is the
template corresponding to the configuration on the right hand side
of Figure \ref{less}

\begin{figure}

\begin{picture}(400,200)
\put(10,0){\begin{picture}(200,200)
\put(80,150){\circle*{6}}		
\put(30,50){\circle*{6}}		
\put(130,50){\circle*{6}}		
\put(80,100){\circle*{6}}		
\put(180,100){\circle*{6}}		

\put(122,35){$x_2$}

\put(10,10){\line(1,2){90}}		
\put(80,-10){\line(0,1){200}}		
\put(150,10){\line(-1,2){90}}		
\put(10,30){\line(1,1){120}}		
\put(10,50){\line(1,0){180}}		
\put(170,10){\line(-1,1){140}}		

\put(200,100){\line(-1,0){180}}		
\put(200,90){\line(-2,1){150}}		
\put(201,107){\line(-3,-1){195}}	
\put(200,120){\line(-1,-1){130}}	
\put(200,96){\line(-5,1){170}}		
\end{picture}}

\put(250,0){\begin{picture}(200,200)
\put(80,150){\circle*{6}}		
\put(30,50){\circle*{6}}		
\put(80,100){\circle*{6}}		
\put(180,100){\circle*{6}}		

\put(10,10){\line(1,2){90}}		
\put(80,-10){\line(0,1){200}}		
\put(10,30){\line(1,1){120}}		

\put(200,100){\line(-1,0){180}}		
\put(200,90){\line(-2,1){150}}		
\put(201,107){\line(-3,-1){195}}	
\end{picture}}

\end{picture}
\caption{}\label{less}
\end{figure}

It is clear that $A_iF^p_q(T) \subseteq F^p_{q-1}(A_i T)$. That
$A_iF^p_q(T)$ is a constructible open subset of $F^p_{q-1}(A_0T)$
follows directly from (\ref{generic}).

The corresponding statements for $E^p_q(T)$ follows from
(\ref{prod}). The statements with $A$ replaced by $B$ follow using
the dual argument.
\end{pf}

\begin{example}
An example where $A_0F^p_q(T)$ is a proper subset of
$F^p_{q-1}(A_0T)$ is given in Figure \ref{face}. If $T$
is the 13-marked template associated to the left hand figure, then
the
right hand configuration is an element of
$F^3_8(A_0T)-A_0(F^3_9(T))$.
\end{example}

\begin{figure}
\begin{picture}(400,200)
\put(10,0){\begin{picture}(200,200)
\put(20,35){\circle*{6}}
\put(20,65){\circle*{6}}
\put(80,35){\circle*{6}}
\put(80,65){\circle*{6}}
\put(35,120){\circle*{6}}
\put(35,180){\circle*{6}}
\put(65,120){\circle*{6}}
\put(65,180){\circle*{6}}
\put(165,35){\circle*{6}}
\put(165,165){\circle*{6}}
\put(150,65){\circle*{6}}
\put(150,135){\circle*{6}}
\put(100,100){\circle*{6}}
\put(80,97){$x_0$}

\put(10,10){\line(1,1){180}}
\put(10,190){\line(1,-1){180}}
\put(150,10){\line(0,1){180}}
\put(10,30){\line(2,1){80}}
\put(10,70){\line(2,-1){80}}
\put(30,110){\line(1,2){40}}
\put(30,190){\line(1,-2){40}}
\end{picture}}

\put(250,0){\begin{picture}(200,200)
\put(20,35){\circle*{6}}
\put(20,65){\circle*{6}}
\put(80,35){\circle*{6}}
\put(80,65){\circle*{6}}
\put(35,120){\circle*{6}}
\put(35,180){\circle*{6}}
\put(65,120){\circle*{6}}
\put(65,180){\circle*{6}}
\put(165,35){\circle*{6}}
\put(165,165){\circle*{6}}
\put(100,65){\circle*{6}}
\put(100,135){\circle*{6}}

\put(10,10){\line(1,1){180}}
\put(10,190){\line(1,-1){180}}
\put(100,10){\line(0,1){180}}
\put(10,30){\line(2,1){80}}
\put(10,70){\line(2,-1){80}}
\put(30,110){\line(1,2){40}}
\put(30,190){\line(1,-2){40}}
\end{picture}}
\end{picture}
\caption{}\label{face}
\end{figure}

\begin{corollary}\label{example}
If $T$ is a $(p+q)$-marked template, then
$$
F^p_q(T) \subseteq F^p_q(A^iA_i T)\quad\text{and}
\quad E^p_q(T) \subseteq E^p_q(A^iA_i T)
$$
for all integers $i$ satisfying $0\le i \le p+q$, and
if $T$ is a $(p+q)$-marked template, then
$$
F^p_q(T) \subseteq F^p_q(B^iB_i T)\quad\text{and}
\quad E^p_q(T) \subseteq E^p_q(B^iB_i T)
$$
for all integers $i$ such that $0\le i \le p+q+1$. \qed
\end{corollary}

We are now ready to prove Theorem \ref{main}.
Throughout the remainder of this section, the ground field
will be $\C$,
unless explicitly stated to the contrary.  Fix $p>0$. The proof
is by induction on $q$.  When $q = 0$, $Y^p_q$ is a point and the
result is trivially true. Now suppose that $q>0$ and that
the result is true for
$Y^p_{q-1}$.  Suppose that $U$ is a non-empty constructible open
subset of $Y^p_q$.  The idea behind the proof is to replace $U$ by
a smaller constructible open set $L$ whose image under $A_0$ is
a constructible open subset of $Y^p_{q-1}$ and such that the map
$L \to A_0L$ is a fibration whose fibers are complements of
arrangements
of hyperplanes in $\P^{p-1}$ of fiber type and whose monodromy
representations are trivial. Using the inductive hypothesis, one then
finds a constructible open subset $L'$ of $A_0L$ which is both a
$K(\pi,1)$ and a rational $K(\pi,1)$. It will then follow from
(\ref{fibration}) that $A_0^{-1}(L')$ is the sought after
constructible open subset of $U$. We now give the details.

Our first task is to find a constructible open subset $W$ of
$Y^p_{q-1}$
such that the restriction of $A_0$ to $U\cap A_0^{-1}W$ is a family
of hyperplane complements where each relative hyperplane is proper
over the base.
There are $(p+q)$-marked templates $T_1,\ldots, T_l$ such that
$$
U = Y^p_q - \bigcup_{j=1}^l F^p_q(T_j).
$$
For each $j$, either $F^p_{q-1}(A_0T_j)$ is all of $Y^p_{q-1}$ or is
a proper closed subvariety. We may suppose that $F^p_{q-1}(A_0T_j)$
is $Y^p_{q-1}$ when $j\le k$ and is a proper subvariety when
$j > k$.
When $j \le k$, set
$$
C_j = Y^p_{q-1} - A_0F^p_q(T_j);
$$
this is a constructible closed proper subset of $Y^p_{q-1}$
by (\ref{generic}). Set
$$
W = Y^p_{q-1} - \left(\bigcup_{j\le k} C_j \cup
\bigcup_{j>k} F^p_{q-1}(A_0 T_j)\right).
$$
Then $W$ is a non-empty constructible open subset of
$Y^p_{q-1}$ and the restriction of  $A_0 : A_0^{-1}W \to W$ to
$F^p_q(T_j)$ is proper and surjective when $j\le k$.

Now
$$
A_0^{-1}W \cap U = A_0^{-1}W - \bigcup_{j\le k} F^p_q(T_j).
$$
is a constructible open subset of $U$.  The fiber of
$$
A_0 : A_0^{-1}W - \bigcup_{j\le k} F^p_q(T_j) \to W
$$
over $(x_1,\ldots,x_{p+q})$ is the complement of an arrangement of
hyperplanes in $\P^{p-1}$ which is derived from the configuration
$(x_1,\ldots,x_{p+q})$ and where each relative hyperplane is proper
over $W$.

Our next task is to replace $W$ by a smaller constructible open
set $O$ such that the restriction of $A_0$ to $U\cap A_0^{-1}O$
is a fiber bundle over $O$.

We say that two linear configurations in $\C^m$ have the {\it same
combinatorics} if their associated partially ordered sets are
isomorphic, or equivalently, if they determine the same template.

\begin{proposition}\label{bundle}
There is a non-empty constructible open subset $O$ of $W$ such that
the restriction of
$$
A_0 : A_0^{-1}W - \bigcup_{j\le k} F^p_q(T_j) \to W
$$
to $A_0^{-1}(O)$ has the property that each of its fibers is the
complement of a linear configuration with the same combinatorics.
Consequently,
$$
A_0 : A^{-1}(O) \to O
$$
is a fiber bundle where the action of $\pi_1(O,\ast)$ on the
homology of the fibers is trivial.
\end{proposition}

\begin{pf}
Let $\F$ be the function field of $Y^p_{q-1}$. Then
$\bigcup_{j\le k} F^p_q(T_j)$ is a configuration of hyperplanes in
$\P^{p-1}(\F)$. For generic $v\in W$, the combinatorics of
the restriction of this configuration to the fiber of $A_0$ over $v$
has the same combinatorics as this configuration over the generic
point of $Y^p_{q-1}$. The set of $v$ for which the combinatorics
is different is a closed constructible subset $F$ of $Y^p_{q-1}$.
The desired constructible open subset of $Y^p_{q-1}$ is then
$O = W-F$.
\end{pf}

Next we further shrink both $U$ and $O$ to make the fibers
of $A_0$ to $U\cap A_0^{-1}O$ complements of arrangements of
hyperplanes of fiber type.

\begin{proposition}
There is a non-empty constructible open subset $O'$ of $O$ and a
$(p+q)$-marked template $T$ such that $A_0F^p_q(T)$ contains $O'$
and such that
$$
F^p_q(T) \supseteq \bigcup_{j\le k} F^p_q(T_j)
$$
and the map
$$
A_0^{-1}O' -  F^p_q(T) \to O'
$$
induced by $A_0$ is a fiber bundle all of whose fibers are
complements of arrangements of hyperplanes of fiber type.
\end{proposition}

\begin{pf}
As in the proof of the previous result, we shall denote the function
field of $Y^p_{q-1}$ by $\F$.  The points $x_1,\ldots,x_{p+q}$ are
defined over $\F$, and therefore may be regarded as a configuration
$\x(\F)$ of points in $\P^{p-1}(\F)$. The set
$$
\bigcup_{j\le k} F^p_q(T_j),
$$
is a configuration of hyperplanes  defined over $\F$ and thus
determines an
element of $\D(\x(\F))$. Let $T'\in \T(\x(\F))$ be the corresponding
template. By (\ref{fiber_type}), there is a template $T\in
\T(\x(\F))$
such that
$$
\H_{T} \supseteq \bigcup_{j\le k} F^p_q(T_j)
$$
and such that $\H_{T}$ is an arrangement of fiber type.  Since $T$ is
defined over the generic point of $Y^p_{q-1}$, it follows from
(\ref{closed}) that $A_0F^p_q(T)$ is a constructible open subset of
$Y^p_{q-1}$. Moreover, the set of $v \in A_0F^p_q(T)$ for which the
combinatorics of the restriction of $\H(T)$ to the fiber of $A_0$
over $v$ is given by $T$ is a constructible open subset of
$A_0F^p_q(T)$. Let $O'$ be the
intersection of this open set with $O$.
\end{pf}

By our inductive hypothesis, the constructible open set
$O'$ of $Y^p_{q-1}$ contains a constructible open
subset $L$ which is a $K(\pi,1)$ and a rational $K(\pi,1)$.
Since $A_0^{-1}L$ is a non-empty constructible open subset of
$Y^p_q$,
$$
V = A_0^{-1}(L) -
\left(F^p_q(T) \cup \bigcup_{j\le k} F^p_q(T_j)\right)
$$
is also a non-empty constructible open subset of $U$. Further,
the map
$$
A_0 : V \to L
$$
is a fibration each of whose fibers is the complement of an
arrangement of hyperplanes in $\P^{p-1}$ of fiber
type. It follows from (\ref{fiber-type}) that the fibers are
$K(\pi,1)$s and rational $K(\pi,1)$s. Since the base is a
$K(\pi,1)$ and a rational $K(\pi,1)$, and since the monodromy is
trivial (\ref{bundle}), it follows from (\ref{fibration})
that $V$ is a $K(\pi,1)$ and a rational $K(\pi,1)$. This
completes the proof of Theorem \ref{main}.

\section{Existence and Uniqueness of Higher Logarithms}
\label{existence}

In this section, we first establish the existence and uniqueness of
the $p$th higher logarithm in the sense of \cite[(6.1)]{hain-macp},
but with $G^p_\bullet$ replaced by a suitably chosen Zariski open
subset $U^p_\bullet$. We then show how to construct the generalized
$p$-logarithm, a multivalued Deligne cohomology class, in the sense
of \cite[(12.4)]{hain-macp}, but with $G^p_\bullet$ replaced by
$U^p_\dot$. We shall use the notation and definitions of
\cite{hain-macp}.

We will say that a simplicial variety $U_\bullet$ is a
{\it subvariety} of the simplicial variety $X_\bullet$ if
each $U_q$ is a subvariety of $X_q$, and if the inclusion $U_\bullet
\hookrightarrow X_\bullet$ is a morphism of simplicial varieties. We
will
say that $U_\bullet$ is an {\it open} (resp.\ {\it closed, dense,
constructible}) subset of $G^p_\bullet$ if each $U_q$ is open (resp.\
{ closed, dense,constructible}) in each $X_q$. There are analogous
definitions with $G^p_\bullet$ replaced by $Y^p_\bullet$.

\begin{proposition}\label{good_sub}
For each positive integer $p$, each dense constructible open subset
$V^p_\bullet$ of the truncated  simplicial variety $G^p_\bullet$
contains
a dense constructible open subset $U^p_\bullet$ where each $U^p_q$
is a rational $K(\pi,1)$. In particular, $G^p_\bullet$ contains a
dense
constructible open subset $U^p_\bullet$ where each $U^p_q$ is a
$K(\pi,1)$ and a rational $K(\pi,1)$.
\end{proposition}

\begin{pf}
The only dense constructible open subset of $G^p_0$ is $G^p_0$
itself.
So $V^p_0=G^p_0$, and we must take $U^p_0=G^p_0$. Suppose that $m >
0$
and that $U^p_q$ has been constructed when $q< m$ such that each
$U^p_q$ is dense in $G^p_q$, $U^p_q \subseteq V^p_q$, and such that
$A_i(U^p_q) \subseteq U^p_{q-1}$ whenever $0 < q < m$.  Now,
$$
V^p_m \cap \bigcap_{i=0}^m A_i^{-1}U^p_{m-1}
$$
is a non-empty constructible open subset of $G^p_m$. So, by
(\ref{main}),
it contains a non-empty, and therefore dense, constructible open
subset
$V^p_m$ of $G^p_m$. The result now follows by induction.
\end{pf}

In order to apply the multivalued de Rham complex functor, we will
need
to know that such a constructible open subset $U^p_\bullet$ of
$G^p_\bullet$
can be lifted to a truncated simplicial object in the category
$\Atilde$
defined in the introduction and in \cite[\S 2]{hain-macp}.

\begin{theorem}\label{lift}
Each constructible open subset $U^p_\bullet$ of $G^p_q$ can be lifted
to a truncated simplicial object of the category $\Atilde$.
\end{theorem}

The lift is natural in the following sense: it comes with a lift
$\tilde{\imath}$ of the inclusion $i:U^p_\dot \hookrightarrow
G^p_\dot$
such that if $j: V^p_\dot\hookrightarrow U^p_\dot$ is an inclusion of
constructible open subsets of $G^p_\dot$, then
$\tilde{\imath}\tilde{\jmath} = \widetilde{\imath\jmath}$.

As the proof of this theorem is technical; it is given in a separate
section, \S \ref{proof}.

Next, we show how to construct the $p$th higher logarithm in
the sense of \cite[(6.1)]{hain-macp} defined on some constructible
dense open subset of $G^p_\bullet$.

The following fact is is a direct consequence of
\cite[(7.8)]{hain-macp}
and \cite[(8.2)(i)]{hain-macp}.

\begin{proposition}\label{acyclic}
If the complex algebraic variety $X$ is a rational $K(\pi,1)$ with
$q(X)=0$, then for all $l\ge 0$, the complex
$W_l\Omegatilde^\bullet(X)$ is acyclic.
\end{proposition}

The existence of the higher logarithms is now an immediate
consequence
of (\ref{good_sub}), (\ref{lift}), (\ref{acyclic}) and
\cite[(9.7)]{hain-macp}:

\begin{theorem}\label{exist_1}
For each integer $p\ge 1$, there is a dense constructible open subset
$U^p_\bullet$ of the simplicial variety $G^p_\bullet$ which has a
lift to the category $\Atilde$, and there is an element $Z_p$ of the
double
complex $W_{2p}\Omegatilde^\bullet (U^p_\bullet)$, unique up to a
coboundary, whose coboundary is the ``volume form''
$$
{dx_1\over x_1} \wedge \ldots \wedge {dx_p\over x_p} \in
\Omega^p(G^p_0).
$$
\end{theorem}

\begin{remark}
With a little more care, one can arrange for each $U^p_q$ to be
invariant under the action of the symmetric group $\Sigma_{p+q+1}$
on $G^p_q$ and for the symbol (as defined in
\cite[p.~444]{hain-macp})
of each component of $Z_p$ to span a copy of the alternating
representation. One should note, however, that it seems difficult to
arrange for each $U^p_q$ to be a rational $K(\pi,1)$ and be
preserved by the action of $\Sigma_{p+q+1}$.
\end{remark}

\begin{remark}\label{hana-mac-const}
Hanamura and MacPherson \cite{hanamura-macp_2} give
an explicit construction of all higher logarithms in the double
complex
$W_{2p}\Omegatilde(G^p_\bullet)$. In particular, they show that it
is not necessary to pass to a Zariski open subset of $G^p_\bullet$
as we did.
\end{remark}

Next, we establish the existence of higher logarithms as Deligne
cohomology classes.  For this, we shall assume the reader is familiar
with the definition of the multivalued Deligne cohomology functor
$\Hmd^\bullet({\blank},\Q(p))$ defined in
\cite[\S 12]{hain-macp}.

The key point here is the following result, a slightly stronger
version
of which was stated in \cite[(12.3)]{hain-macp}, and which we will
prove in \S \ref{pf_dunno}. Recall from the introduction that
$\Hdel^\bullet$ denotes Beilinson's absolute Hodge cohomology.

\begin{theorem}\label{dunno}
Suppose that $X_\bullet$ is a truncated simplicial variety with
a lift to $\Atilde$. If each $X_q$ is a rational $K(\pi,1)$, then
for each integer $p$, there is a natural isomorphism
$$
\Hmd^\bullet(X_\bullet,\Q(p)) \approx
\Hdel^\bullet(X_\bullet,\Q(p))).
$$
\end{theorem}

Granted this, the construction of the generalized $p$th higher
logarithm as an element of $\Hmd^\bullet(U^p_\bullet,\Q(p))$ is
relatively straightforward.

\begin{theorem}\label{gen_log}
If $U^p_\bullet$ is a dense subvariety of $G^p_\bullet$ where each
$U^p_q$ is a rational $K(\pi,1)$, then there is an element of
$$
\Hmd^{2p}(U^p_\bullet,\Q(p))
$$
whose restriction to $G^p_0$ is the volume form.
\end{theorem}

\begin{pf}
Let $V^p_m$, be the subvariety of
$$
\left\{(v_0,v_1,\ldots,v_{m}): v_j \in \C^p\right\}
$$
which consists of those $(m+1))$-tuples of vectors where each set of
$\min(m+1,p)$ of the vectors is linearly independent. When
$m\ge p$, there is a natural projection $V^p_m\to G^p_{m-p}$ which
is a principal $GL_p(\C)$-bundle. Define face maps
$$
A_i : V^p_m \to V^p_{m-1}
$$
by omitting the $i$th vector. Denote the corresponding simplicial
variety by $V^p_\bullet$. Denote the truncated simplicial space
which consists only of those $V^p_m$ with $p\le m \le 2p$ by
$\Vtilde^p_\bullet$. There is a natural projection
$\Vtilde^p_\bullet \to G^p_\bullet$ which is a principal
$GL_p(\C)$-bundle. We would like to say that this bundle has a
Chern class
$$
c_p \in \Hdel^{2p}(G^p_\bullet,\Q(p)).
$$
Since the variety $G^p_\bullet$ is truncated (it has no simplices
in dimensions $< p$), the existence of such a Chern class is not
immediate. Our next task is to establish the existence of this class.
We do this using the Borel construction.

Let $E_\bullet$ be, say, the standard simplicial model for the
universal bundle associated to $GL_p(\C)$. What is important for us
is that $E_\bullet$ is a simplicial variety with the homotopy type
of a point and on which $GL_p(\C)$ acts freely. Let $P_\bullet$
be the bisimplicial variety $V_\bullet \times E_\bullet$.  It
has the homotopy type of $V^p_\bullet$.  Denote the quotient of
$P_\bullet$ by the diagonal action of $GL_p(\C)$ by $B_\bullet$.
Since $GL_p(\C)$ acts freely on $P_\bullet$, the quotient map
$$
P_\bullet \to B_\bullet
$$
is a principal $GL_p(\C)$ bundle. By \cite{beilinson}, this bundle
has a Chern class
\begin{equation}\label{chern}
c_p \in \Hdel^{2p}(B_\bullet,\Q(p)).
\end{equation}

Denote the truncated simplicial variety consisting of those $B_m$
with $p\le m \le 2p$ by $\Btilde_\bullet$ and denote the restriction
of the bundle $P_\bullet$ to $\Btilde_\bullet$ by $\Ptilde_\bullet$.
Then there is a commutative diagram
$$
\begin{matrix}
\Ptilde_\bullet & \to & \Vtilde^p_\dot \\
\downarrow & & \downarrow \\
\Btilde_\dot & \to & G^p_\bullet\\
\end{matrix}
$$
of principal $GL_p(\C)$ bundles obtained by collapsing out
$E_\bullet$. Since the action of $GL_p(\C)$ on $\Vtilde^p$ is free,
the bottom arrow is a homotopy equivalence of simplicial varieties,
and therefore induces an isomorphism on Deligne cohomology.

We can therefore restrict the Chern class (\ref{chern}) to
$G^p_\bullet$ to obtain a class in
$$
\Hdel^{2p}(G^p_\bullet,\Q(p)).
$$
It follows from (\ref{dunno}) that we can restrict this class to
$U^p_\bullet$ to obtain a class $C_p$ in
$$
\Hmd^{2p}(U^p_\bullet,\Q(p))
$$
provided that each $U^p_q$ is a rational $K(\pi,1)$.

Finally, to prove Theorem \ref{gen_log}, we have to show that the
restriction of $C_p$ to $U^p_0=G^p_0$ is a non-zero multiple of
the volume form in $H^p(G^p_0)$. It is proved in \cite{hain-yang}
that the restriction of $C_p$ to $G^p_0$ is $(p-1)! \vol$. It follows
that $C_p/(p-1)!$ is a generalized $p$-logarithm. This completes
the proof of Theorem \ref{gen_log}.
\end{pf}

\section{Proof of Theorem \ref{lift}}
\label{proof}

In the proof we shall need the following construction. Let
$$
X^p_q = \C^{p+q+1} - \Delta,
$$
where $\Delta$ denotes the fat diagonal---that is, the locus of
points
in $\C^{p+q+1}$ where the coordinates are not all distinct. Define
$$
A_i : X^p_q \to X^p_{q-1}
$$
by deleting the $i$th coordinate:
$$
A_i: (t_0,\ldots,t_{p+q}) \mapsto
(t_0,\ldots,\widehat{t_i},\ldots,t_{p+q}).
$$

Denote the truncated simplicial variety consisting of those $X^p_q$
with
$0\le q \le p$ by $X^p_\bullet$. We can define a morphism
$\phi : X^p_\bullet \to G^p_\bullet$ by taking $(t_0,\ldots,t_{p+q})$
to the $GL_p(\C)$ orbit of the $(p+q+1)$-tuple of vectors
$$
\begin{pmatrix}
1 \cr t_{0} \cr  t_{0}^2 \cr \vdots \cr t_{0}^{p-1}
\end{pmatrix}
,\quad
\begin{pmatrix}
1 \cr t_{1} \cr  t_{1}^2 \cr \vdots \cr t_{1}^{p-1}
\end{pmatrix}
,\,\dots \quad
\begin{pmatrix}
1 \cr t_{p+q} \cr  t_{p+q}^2 \cr \vdots \cr t_{p+q}^{p-1}
\end{pmatrix}.
$$
This map is easily seen to be a well defined morphism of simplicial
varieties. It induces a morphism $\phibar : X^p_\bullet \to
Y^p_\bullet$.

\begin{lemma}\label{dense}
The image of $X^p_q$ in $G^p_q$ is dense in $G^p_q$ in the
constructible topology.
\end{lemma}

\begin{pf}
In view of (\ref{prod}), we need only prove that the image of
$X^p_q$ in $Y^p_q$ is dense in $Y^p_q$ in the constructible topology.
We do this by induction on $q$. Since $Y^p_0$ is a point, the result
is trivially true when $q=0$. Suppose that $q>0$. Denote the
constructible
closure of the image of $X^p_q$ in $Y^p_q$ by $C^p_q$. By induction,
$C^p_{q-1} = Y^p_{q-1}$. The intersection of $C^p_q$ with each fiber
of $A_0 : Y^p_q \to Y^p_{q-1}$ is a constructible closed subset of
the
fiber. Note that the fiber of $A_0 : Y^p_q \to Y^p_{q-1}$ is the
complement of a linear arrangement in $\P^{p-1}$ and that each of its
constructible closed subsets is the intersection of the fiber with a
finite union of linear subspaces of $\P^{p-1}$. Since
the intersection of $C^p_q$ with the fiber is an open subset of a
rational normal curve in $\P^{p-1}$, and since each rational normal
curve is non-degenerate, it follows that the fiber of
$A_0 : C^p_q \to Y^p_{q-1}$ equals the fiber of
$A_0 : Y^p_q \to Y^p_{q-1}$. It follows that $C^p_q = Y^p_q$.
\end{pf}

Denote the topological analogue of the category $\Atilde$ by
$\toptilde$.
Observe that a simplicial object of $\A$ has a lift to the category
$\Atilde$ if and only if it has a lift to the category $\toptilde$.

\begin{proposition}
Suppose that $Y_\dot$ and $Z_\dot$ are simplicial
topological spaces where each $Y_n$ and $Z_n$ is path connected.
If $f : Y_\dot \to Z_\dot$ is a morphism of simplicial spaces, and if
$Y_\dot$ has a lift to $\toptilde$, then both $Z_\dot$ and $f$ have
lifts to $\toptilde$.
\end{proposition}

\begin{pf}
We use the equivalence of the category $\Atilde$ with the category
$\A_\ast$ which is constructed in \cite[\S 2]{hain-macp}. We first
construct a simplicial object of $\A_\ast$ which corresponds to
the lift of $Y_\dot$ to a simplicial object of $\Atilde$.

Let $\Ytilde_\dot$ be the simplicial object of $\Atilde$ which
is the lift of $Y_\dot$.  Choose a base point $y_n'$ of $\Ytilde_n$
for each
$n$, and let $y_n$ be its image in $Y_n$.  Each strictly order
preserving map $\phi : [m] \to [n]$ induces a morphism
$A_\phi' : \Ytilde_n \to \Ytilde_m$ of $\Atilde$ which covers the
face map $A_\phi : Y_n \to Y_m$. Since each
$\Ytilde_n$ is connected and simply  connected, there is a unique
homotopy class of paths in $\Ytilde_m$ from $y_m'$ to
$A_\phi'(y_n')$. Its image in $Y_m$ is a distinguished homotopy
class of paths $\gamma_\phi$ in $Y_m$ from $y_m$ to
$A_\phi(y_n)$. The pair $(A_\phi,\gamma_\phi)$ is a morphism
$(Y_n,y_n) \to (Y_m,y_m)$ in the category $\A_\ast$ and the
collection $(Y_n,y_n)$ of pointed spaces together with the maps
$(A_\phi,\gamma_\phi)$ is a simplicial object of $\A_\ast$.

We now use this to construct a lift of $X_\dot$ to $\A_\ast$. Let
$x_n = f(y_n)$. For each order preserving injection $\phi : [m] \to
[n]$, let $\mu_\phi$ be the homotopy class $f\cdot \gamma_\phi$ of
paths
in $X_m$ from $x_m$ to $A_\phi(x_n)$. The collection of pointed
spaces
$(X_n,x_n)$ together with the pairs $(A_\phi,\mu_\phi)$ is easily
seen to be a simplicial object of $\A_\ast$.  Take $\Xtilde_n$ to be
the standard model of the universal covering space of $(X_n,x_n)$
--- it consists of homotopy classes $\rho$ of paths that emanate from
$x_n$. The face maps $A_\phi$ lift to face maps $A_\phi'$ by defining
$A_\phi'(\rho)$ to be the homotopy class of paths $\mu_\phi\cdot
\rho$
in $\Xtilde_m$. This is a simplicial object of $\Atilde$ which lifts
$X_\dot$.
\end{pf}

\begin{corollary}\label{top-lift}
Suppose that $Y_\dot$ and $Z_\dot$ are simplicial
topological spaces where each $Y_n$ and $Z_n$ is path connected.
If $f : Y_\dot \to Z_\dot$ is a morphism of simplicial spaces, and if
each simplex $Y_n$ of $Y_\dot$ is simply connected, then $Z_\dot$ has
a canonical lift to $\toptilde$ such that $f$ is a morphism of
$\toptilde$.\qed
\end{corollary}

The following result is needed in the proof of the theorem.

\begin{lemma}\label{poly-lem}
Suppose that $f \in \R[t_1,\dots,t_n]$. If $f\neq 0$, there is a
real number $K > 1$ such that $f$ is bounded away from zero in the
region
$$
D_n(K) := \left\{(t_1,\dots,t_n) :
t_1 \ge K, t_2 \ge Ke^{t_1}, \dots t_n \ge Ke^{t_{n-1}}\right\}.
$$
\end{lemma}

\begin{pf}
The proof is by induction on the number of variables $n$. The result
is
trivially true when $n=1$. Now suppose that $n>1$ and that the result
has been proved for polynomials with fewer than $n$ variables. Set
$x = (t_1,\dots, t_{n-1})$ and $y = t_n$. We can write
\begin{equation}\label{poly}
f = a_d(x)y^d + a_{d-1}(x)y^{d-1} + \dots + a_1(x)y + a_0(x)
\end{equation}
where each $a_j(x) \in \R[t_1,\dots,t_{n-1}]$ and $a_d\neq 0$. If
$d=0$,
then we are in the previous case and the result holds by induction.
So
assume that $d>0$.  By induction, there exist real constants $C>0$
and
$K>1$ such that $ |a_d(x)| \ge C $ for all $x\in D_{n-1}(K)$. By a
standard estimate, the roots $\theta(x)$ of the polynomial
(\ref{poly}) satisfy
$$
|\theta(x)| \le 1 + \max_{0 \le j < d} \left|{a_j(x) \over
a_d(x)}\right|
\le 1 + \max_{0 \le j < d} {|a_j(x)| \over C} \le \|x\|^l
$$
for some positive integer $l$ and for each $x \in D_{n-1}(K)$.

Observe that if $(t_1,\dots,t_{n-1}) \in D_{n-1}(K)$, then
$$
1 < K \le t_1 \le t_2 \le \dots \le t_{n-1}
$$
so that $\|x\| \le t_{n-1}$ when $x\in D_{n-1}(K)$.
It follows from the previous inequality that
$$
1 + |\theta(x)| < e^{t_{n-1}} \le y
$$
provided that $t_{n-1}$ is sufficiently large, which can be arranged
by
increasing $K$ if necessary. Since
$$
f(x,y) = \pm\, a_d(x) \prod_{j=1}^d (y-\theta_j(x)),
$$
it follows that if $(x,y) \in D_n(K)$, then $|f(x,y)| \ge C$.
\end{pf}

\begin{pf*}{Proof of Theorem \ref{lift}}
We first give a brief proof of (\ref{lift}) in the case when
$U^p_\dot = G^p_\dot$. For each $q$, the subset
$$
\Delta^p_q := \left\{(t_0,\ldots,t_{p+q}) :
t_i \in \R, 0 < t_0 < t_1 < \cdots < t_{p+q}\right\}
$$
of $X^p_q(\R)$ is contractible. Moreover, each of the face maps
$A_i$ maps $\Delta^p_q$ into $\Delta^p_{q-1}$. It follows that we
have morphisms
$$
\Delta^p_\bullet \hookrightarrow X^p_\bullet \to G^p_\bullet
$$
of truncated simplicial spaces. Since each $\Delta^p_q$ is
contractible,
$\Delta^p_\bullet$ has a unique lift to a simplicial object of the
category $\toptilde$. If follows from (\ref{top-lift}) that both
$X^p_\dot$ and $G^p_\dot$ have lifts to $\toptilde$, and therefore to
$\Atilde$.

The strategy in the general case is similar. We seek a simplicial
space
$D_\dot$, each of whose simplices is contractible, which maps to
$U^p_\dot$. It follows from (\ref{dense}) that the pullback of
$U^p_q$
to $X^p_q$ is a proper open subvariety$V^p_q$ of $X^p_q$. By standard
arguments, there is a non zero polynomial $f_q\in \R[t_0,\dots,
t_{p+q}]$
such that
$$
X^p_q - f_q^{-1}(0) \subseteq V^p_q.
$$
It follows from (\ref{poly-lem}) that there is a real number
$K_q > 1$ such that
$$
D^p_q(K_q):=
\left\{(t_0,\dots,t_{p+q}) \in X^p_q(\R) : t_0 \ge K_q, t_j \ge
K_qe^{t_{j-1}} \text{ when } j\ge 1 \right\} \subset V^p_q.
$$
Let
$$
K = \max_{0\le q \le p} K_q.
$$
Then $D^p_q(K) \subset V^p_q$.
It is not difficult to show that $A_i(D^p_q(K)) \subseteq D^p_{q-1}$
for
each $i$. It follows that the $D^p_q(K)$, with $0\le q \le p$, form a
truncated simplicial space $D^p_\dot(K)$ which maps to $V^p_\dot$,
and
therefore to $U^p_\dot$. It is not difficult to show that each
$D^p_q(K)$
is contractible. It follows from (\ref{top-lift}) that $U^p_q$ lifts
to
a simplicial object of $\Atilde$.
\end{pf*}

\section{Proof of \ref{dunno}\label{pf_dunno}}

We only give a detailed sketch of the proof. First we prove the
result
when $X_\dot$ is replaced by a single space.

Denote the Malcev Lie algebra associated to the pointed space
$(Y,y)$ by $\p(Y,y)$. Now suppose that $Y$ is a complex algebraic
manifold. Recall from \cite[p.~470]{hain-macp} that
the multivalued Deligne cohomology of the object $(Y,y)$ of the
category $\Atilde$ is defined to be the cohomology of the
complex\footnote{Note that there is a typo in the definition of
$MD(X,\Q(p))$ in \cite[p.~470]{hain-macp} --- one should quotient
out by $F^p\calC(\g,\Omegatilde^\dot)$ as defined on op cit, p.~469
and not just by $F^p\Omegatilde^\dot$.}
\begin{multline*}
MD(Y,\Q(p)) :=  \cone\Bigl(
W_{2p} \Homcts_\Q(\Lambda^\dot \p(Y,y),\Q) \oplus \\
W_{2p}\Homcts_\C(\Lambda^\dot \p(Y,y),F^p\Omegatilde^\dot(Y,y))
\to W_{2p}\Homcts_\C(\Lambda^\dot \p(Y,y),\Omegatilde^\dot(Y,y))
\Bigr)[-1].
\end{multline*}

We shall need a $\Q$ analogue of
$\Homcts_\C(\Lambda^\dot \p(Y,y),\Otilde(Y))$.
This will be constructed using continuous cohomology of certain path
spaces.

The space of paths in a topological space $Y$ which go from $y\in Y$
to $z\in Y$ will be denoted by $P_{y,z}Y$. The homology group
$H_0(P_{y,z}Y,\Q)$ has a natural topology which agrees with the
filtration of $H_0(P_{y,y}Y,\Q)\approx\Q\pi_1(Y,y)$ by powers of
its augmentation ideal when
$y=z$ (cf. \cite[\S 3]{hain-zucker}). Denote the continuous dual
of $H_0(P_{y,z}Y,\Q)$ by
$\Hcts^0(P_{y,z}Y,\Q)$.  These groups fit together to form a local
system over $Y\times Y$ whose fiber over $(y,z)$ is
$\Hcts^0(P_{y,z}Y,\Q)$. It is a direct limit of unipotent local
systems
over $Y\times Y$ and a direct limit of unipotent variations of
mixed Hodge structure when $Y$ is a smooth algebraic variety
\cite{hain-zucker}.

There are two  natural inclusions of
$\Homcts_\Q(\Lambda^\dot \p(Y,a),\Q)$, ($a=y,z$),
into
$$
\Homcts_\Q
(\Lambda^\dot \p(Y,y) \otimes_\Q \Lambda^\dot
\p(Y,z),\Hcts^0(P_{y,z}Y)).
$$
They are induced by the two projections of $(Y,y)\times (Y,z)$ onto
$(Y,a)$ and by the inclusion of the constants into
$\Hcts^0(P_{y,z}Y)$.
We shall denote them by $\phi_1$ and $\phi_2$, respectively.

\begin{proposition}
If $H_1(Y,\Q)$ is finite dimensional, then $\phi_1$ and $\phi_2$
are both quasi-isomorphisms.
\end{proposition}

\begin{pf}
We prove the result for $\phi_2$, the other case being similar.
By a standard spectral sequence argument, it suffices to show that
$$
\Homcts_\Q(\Lambda^\dot \p(Y,y), \Hcts^0(P_{y,z}Y))
$$
is acyclic. We may, without loss of generality, take $y=z$.
Because $H_1(Y)$ is finite dimensional, each graded quotient of the
topology on $H_0(P_{y,y}Y,\Q)$ is finite dimensional, and it follows
that the dual of $\Hcts^0(P_{y,y}Y),\Q)$ is isomorphic to the
completion
of $H_0(P_{y,y}Y,\Q)\approx \Q\pi_1(Y,y)$. This, in turn, is
isomorphic
to the completion of $U\p(Y,y)$.
So there is a natural isomorphism
$$
\Homcts_\Q(\Lambda^\dot\p(Y,y),\Hcts^0(P_{y,y}Y))
\approx \Homcts_\Q(\Lambda^\dot\p(Y,y) \comptensor U\p(Y,y),\Q)
$$
of chain complexes. This last complex is acyclic, as it is the
continuous dual of an acyclic complex (cf.\
\cite[(3.9)]{hain:cycles}).
\end{pf}

The following result is a straightforward refinement of the previous
result.

\begin{proposition}\label{qism}
If $Y$ is a smooth algebraic variety, then each of the complexes in
the previous result is a a complex of mixed Hodge structures, and the
two natural inclusions
$\phi_1$ and $\phi_2$ of $\Homcts_\Q(\Lambda^\dot \p(Y,a),\Q)$
($a=y,z$) into
$$
\Homcts_\Q(\Lambda^\dot \p(Y,y) \otimes_\Q \Lambda^\dot \p(Y,z),
\Hcts^0(P_{y,z}Y))
$$
are quasi-isomorphisms in the category of complexes of mixed Hodge
structures. \qed
\end{proposition}

Next, observe that each $F \in\Otilde(Y,y)$ induces a linear map
$$
H_0(P_{y,z}Y) \to \C
$$
by taking the path $\gamma$ to the difference $F(z) - F(y)$ where the
branch of $F$ at $z$ is obtained by analytically continuing $F$ along
$\gamma$. It follows from standard properties of iterated integrals
that this map is continuous. Consequently, we obtain a linear
map
$$
\Otilde(Y,y) \to \Hcts^0(P_{y,z}Y,\C).
$$
It follows from \cite[\S 3]{hain-macp} and Chen's de Rham Theorem
for the fundamental group that when $q(Y)=0$, this map is an
isomorphism of  $W_\dot$ filtered vector spaces. This isomorphism
is $\pi_1(Y,y)$-equivariant with respect to the standard actions
of $\pi_1(Y,y)$ on $\Otilde(Y,y)$ and $\Hcts^0(P_{y,z}Y,\C)$.

Recall that there is a natural homomorphism
$$
\theta : \Homcts_\C(\Lambda^\dot \p(Y,y), \C) \to \Omega^\dot(Y)
$$
of $W_\dot$ filtered d.g.\ algebras \cite[(7.7)]{hain-macp}.

Fix a point $z$ of $Y$. Consider the complex
\begin{multline*}
\cone\Bigl(W_{2p}
\Homcts_\Q\left(\Lambda^\dot \p(Y,y) \otimes_\Q \Lambda^\dot \p(Y,z),
\Hcts^0(P_{y,z}Y,\Q)\right) \oplus \\
W_{2p}\Homcts_\C(\Lambda^\dot \p(Y,y),F^p\Omegatilde^\dot(Y,y))
\to W_{2p}\Homcts_\C(\Lambda^\dot \p(Y,y),\Omegatilde^\dot(Y,y))
\Bigr)[-1].
\end{multline*}
Here
$$
\Homcts_\Q(\Lambda^\dot \p(Y,y) \otimes_\Q \Lambda^\dot \p(Y,y),
\Hcts^0(P_{y,z}Y,\Q))
$$
is mapped into
$$
W_{2p}\Homcts_\C(\Lambda^\dot \p(Y,y),\Otilde(Y,y)\otimes
\Omega^\dot(Y))
$$
using the identification of $\Hcts^0(P_{y,z}Y,\C))$ with
$\Otilde(Y,y)$
and the map $\theta$ in the second factor. It is straightforward to
check
it is a chain map.

Define a map from this complex to $MD(Y,\Q(p))$ by defining it to be
$\phi_1$ on the first factor and the identity on the other two
factors.
Since $\phi_1$ is a $W_\dot$ filtered quasi-isomorphism, this map is
a
quasi-isomorphism.

Next, Define $MD'(X,\Q(p))$ to be the complex
$$
\cone\Bigl(W_{2p}\Homcts_\Q(\Lambda^\dot \p(Y,y),\Q)
\oplus F^p W_{2p}\Omega^\dot(Y))
\to W_{2p}\Omega^\dot(Y)) \Bigr)[-1],
$$
where the map $\Homcts_\Q(\Lambda^\dot \p(Y,y),\Q) \to
\Omega^\dot(Y)$
is induced by $\theta$.
It can be mapped to the previous complex using $\phi_2$ on the first
factor and the obvious inclusions on the other two factors. Since
$$
\Homcts_\C(\Lambda^\dot \p(Y,y),\Otilde^\dot(Y,y))
$$
is an acyclic complex of mixed Hodge structures, it follows that this
map is also a quasi-isomorphism.
That is, we can equally well compute $\Hmd^\dot(Y,\Q(p))$ using
the complex $MD'(Y,\Q(p))$.

A map of $MD'(Y,\Q(p))$ into a standard complex that computes
$\Hdel^\dot(Y,\Q(p))$ can now be constructed using the techniques of
the proof of \cite[(11.7)]{carlson-hain}. Taking homology, we obtain
a map
$$
\psi : \Hmd^\dot(Y,\Q(p)) \to \Hdel^\dot(Y,\Q(p))
$$
for all smooth varieties.

\begin{proposition}\label{iso}
If $Y$ is a rational $K(\pi,1)$, then $\psi$ is an isomorphism.
\end{proposition}

\begin{pf}
The homology of the complex $\Homcts_\Q(\Lambda^\dot \p(Y,y),\Q)$ is
the continuous cohomology $\Hcts^\dot(\p(Y,y))$ of the Lie algebra
$\p(Y,y)$. The natural map
$$
\Hcts^\dot(\p(Y,y)) \to H^\dot(Y,\Q),
$$
is a morphism of mixed Hodge structures \cite[(11.7)]{carlson-hain}.
Since the multivalued Deligne cohomology is constructed as a cone,
we have a long exact sequence
$$
 \cdots \to W_{2p}\Hcts^{k-1}(Y,\Q) \oplus F^pW_{2p}\Omega^{k-1}(Y)
\stackrel{\theta - i}{\longrightarrow} W_{2p}\Omega^{k-1} (Y)
\to \Hmd^k(Y,\Q(p))
\to \cdots
$$
where $i$ denotes the inclusion of $F^p\Omega^\dot$ into
$\Omega^\dot$.
The map $\psi$ induces a map from this long exact sequence into the
standard long exact sequence
$$
 \cdots \to W_{2p}H^{k-1}(Y,\Q) \oplus F^pW_{2p}H^{k-1}(Y)
\to W_{2p}H^{k-1}(Y,\C) \to \Hdel^k(Y,\Q(p)) \to \cdots
$$
When $Y$ is a rational $K(\pi,1)$, each of the maps
$$
\Hcts^\dot(Y,\Q) \stackrel{\theta}{\to} \Omega^\dot(Y) \to
H^\dot(Y,\C)
$$
is a $(W_\dot,F^\dot)$ bifiltered quasi-isomorphism
\cite[(8.2)(i),(iii)]{hain-macp}. The result now follows using the
5-lemma.
\end{pf}

One can take $y=z$ in each of the chain maps above. If one does this,
the assignment of each of these complexes to an object of the
category
$\A_\ast$ defined in \cite[\S 2]{hain-macp} is a functor.

Suppose that $X_\dot$ is a simplicial object of the category
$\Atilde$.
Choose a base point $x_n$ of each $X_n$; $X_\dot$ now determines a
simplicial object of the category $\A_\ast$, and we may apply any of
the functors above to $X_\dot$ to obtain a double complex. Using
standard
arguments, we see that the total complex associated to the double
complex
$MD'(X_\dot,\Q(p))$ computed the multivalued Deligne cohomology of
$X_\dot$
and that there is a map
$$
\Psi : \Hmd^\dot(X_\dot,\Q(p)) \to \Hdel^\dot(X_\dot,\Q(p)).
$$
When each $X_n$ is a rational $K(\pi,1)$, it is not difficult to
show, using an argument similar to the proof of (\ref{iso}) and the
skeleton filtration, that $\Psi$ is an isomorphism.

\section{Higher Logarithms and Extensions of Tate Variations}
\label{conclusion}

The second remark is that the higher logarithms we have constructed
generically on $G^p_\dot$ are related to extensions of (Tate)
variations
of mixed Hodge structures. Indeed, by \cite[(12.1)]{carlson-hain} and
\cite[(8.6)]{hain:cycles}, if a space $X$ is a rational $K(\pi,1)$
with
$q(X)=0$, then there are natural isomorphisms
$$
\Hdel^\dot(X,\Q(p)) \approx \Ext^\dot_{\H(X)}(\Q,\Q(p)) \approx
\Ext^\dot_{\T(X)}(\Q,\Q(p))
$$
where $\H(X)$ and $\T(X)$ denote the categories of unipotent
variations of mixed Hodge structure over $X$ and Tate variations of
mixed Hodge structure over $X$, respectively. Thus if $X_\dot$ is a
simplicial variety where each $X_n$ is a rational $K(\pi,1)$ with
$q(X)=0$, then, in some sense, we may identify
$\Hmd^\dot(X_\dot,\Q(p))$
with the ``hyper-ext" group of extensions of $\Q$ by $\Q(p)$
associated
to $X_\dot$ (cf. \cite[\S 10]{hain:cycles}).

\end{document}